\documentclass[letterpaper,11pt]{article}

\usepackage[thmmarks]{ntheorem}   
\usepackage{amsmath}              
\usepackage{mathrsfs}		
\usepackage{accents}              
\usepackage{amsopn,amssymb}       
\usepackage{bbm}                  
\usepackage{enumitem}            
\usepackage[margin=1in]{geometry} 
\usepackage{microtype}            
\usepackage{multicol}             
\usepackage{shuffle}              
\usepackage{soul}
\usepackage{tikz-cd}              
\usepackage{verbatim}
\usepackage{titlesec}             

\DeclareMathOperator{\Gal}{Gal}

\DeclareMathOperator{\rig}{rig}

\newcommand{\bC}{\mathbb{C}}

\newcommand{\bQ}{\mathbb{Q}}
\newcommand{\bR}{\mathbb{R}}

\newcommand{\bZ}{\mathbb{Z}}

\titleformat{\section}{\normalfont\scshape\filcenter}{\thesection}{.2cm}{}[\vspace{-.3cm}]
\titleformat{\subsection}[block]{\normalfont\scshape}{\thesubsection}{.2cm}{}

\theoremheaderfont{\kern-0cm\normalfont\bfseries}\theorembodyfont{\slshape}

\theoremstyle{plain}
\theoremseparator{.}

\newtheorem{Theorem}{Theorem}[section]

\newtheorem{Lemma}[Theorem]{Lemma}
\newtheorem{Proposition}[Theorem]{Proposition}

\newtheorem{Definition/Theorem}[Theorem]{Definition/Theorem}
\newtheorem{Definition/Proposition}[Theorem]{Definition/Proposition}

\theorembodyfont{\upshape}
\theoremseparator{.}
\newtheorem{Example}[Theorem]{Example}
\newtheorem{Remark}[Theorem]{Remark}
\newtheorem{Algorithm}[Theorem]{Algorithm}

\theoremstyle{nonumberplain}

\theoremstyle{emptybreak}
\theoremheaderfont{\kern-0cm\normalfont\bfseries}\theorembodyfont{\slshape}

\theoremstyle{plain}
\theoremseparator{.}
\newtheorem{Definition}[Theorem]{Definition}

\theoremstyle{nonumberplain}
\theoremheaderfont{\kern-0cm\sc}\theorembodyfont{\upshape}
\theoremseparator{}
\theoremsymbol{\ensuremath{\square}}
\newtheorem{Proof}{Proof.}
\newtheorem{nProof}{}
\usepackage{arydshln}
\usepackage{xcolor}
\usepackage{graphicx}
\usepackage{hyperref}
\usepackage{float}
\usepackage[font = footnotesize]{caption}

\title{Searching for Rigidity in Algebraic Starscapes}
\author{Gabriel Dorfsman-Hopkins, Shuchang Xu}
\date{}
\begin{document}
\maketitle
\begin{abstract}
  We create plots of algebraic integers in the complex plane, exploring the effect of sizing the points according to various arithmetic invariants.  We focus on Galois theoretic invariants, in particular creating plots which emphasize algebraic integers whose Galois group is \emph{not} the full symmetric group$-$these integers we call \emph{rigid}.  We then give some analysis of the resulting images, suggesting avenues for future research about the geometry of so-called \emph{rigid algebraic integers}.
\end{abstract}
\begin{figure}[h!]
  \centering
  \includegraphics[scale=.2]{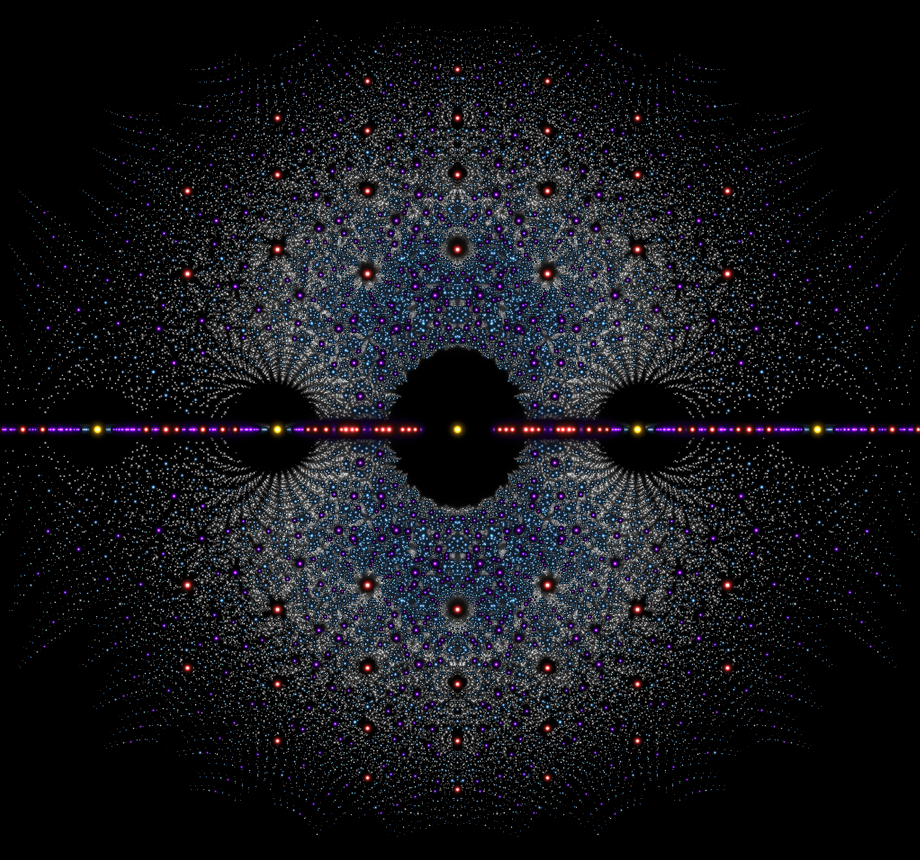}\hspace{20pt}\includegraphics[scale=.2]{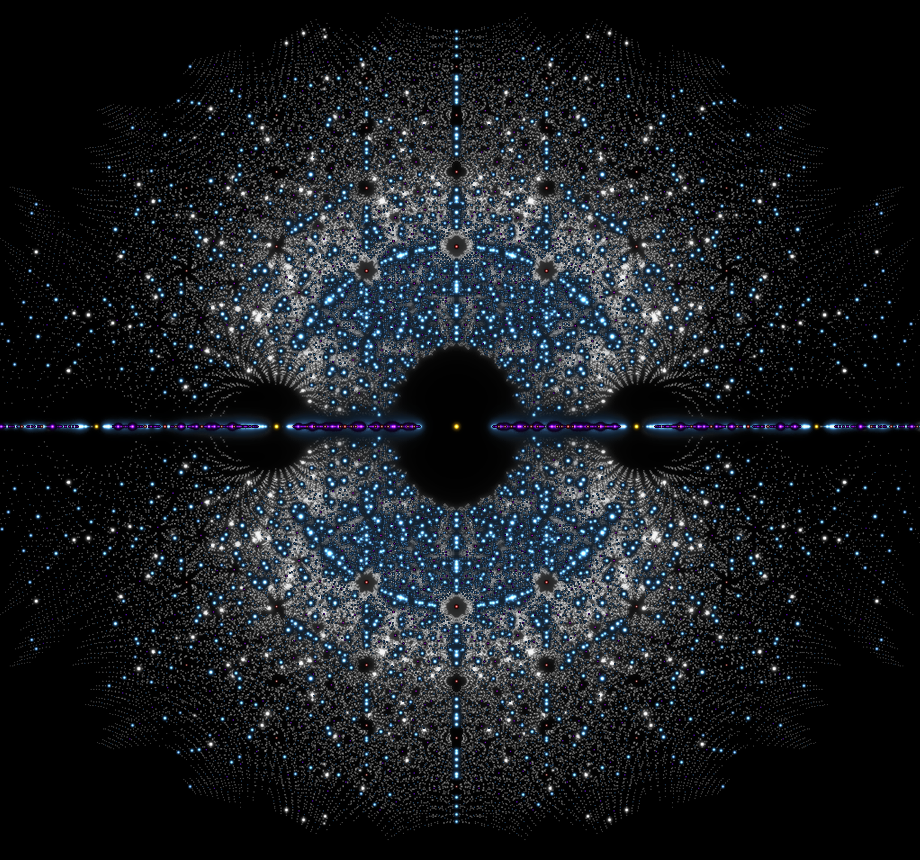}\\
  \vspace{20pt}
  \includegraphics[scale=.2]{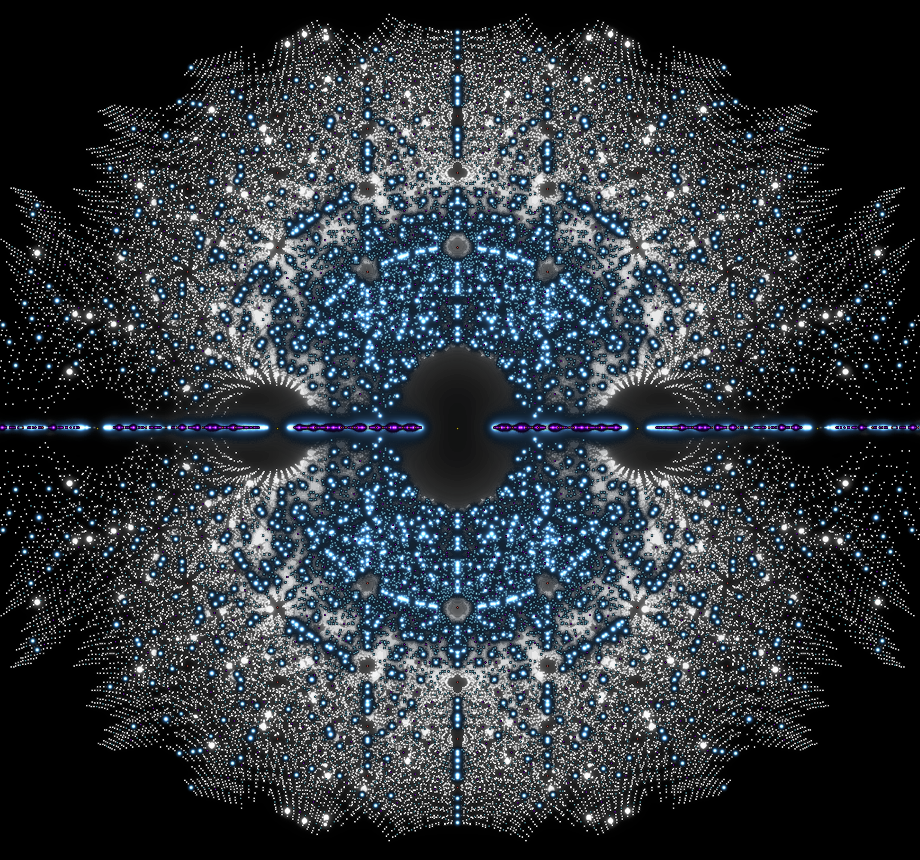}\hspace{20pt}\includegraphics[scale=.2]{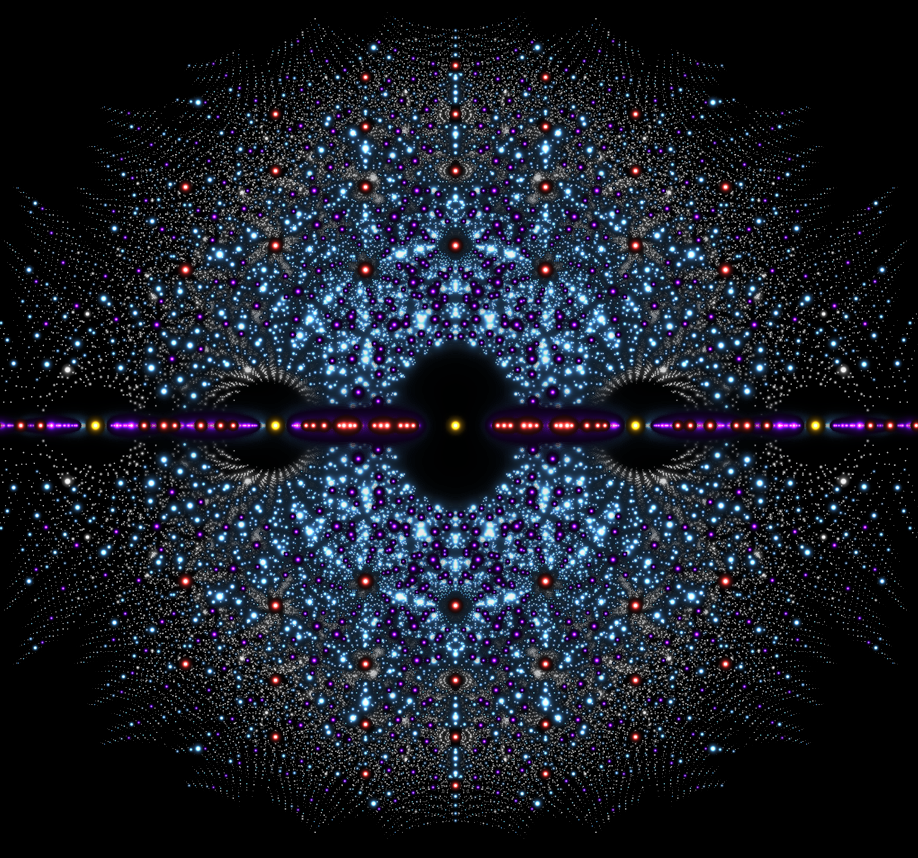}
  \caption{Perspectives on the algebraic integers of degree $\le5$.  Brightness parameters (as defined in Section \ref{imageGeneration}) are as follows.  \textit{Top left:} Sized by discriminant only, using $\beta_D$ as the brightness parameter.  \textit{Top right:} by coefficient vector and rigidity, using $\beta_{vr}$.  \textit{Bottom left:} by rigidity only, using $\beta_r$.  \textit{Bottom right}: by discriminant and rigidity, using $\beta_{Dr}$.}
  \label{Frontpage}
\end{figure}

\section{Introduction}\label{intro}
In what follows, we visually explore the set of complex numbers that are roots of monic polynomials with integer coefficients, the so-called \emph{algebraic integers}.  The creation and study of plots of algebraic integers has a rich and collaborative history, bringing together pure and computational mathematics with digital art.  These images exhibit deep relationships between geometry and arithmetic, and serve as invitations to explore the mysterious patterns lying within the integers, while attempts to explain the visual structure in these plots have inspired deep research programs.

Dating back to the early 90s, Odlyzko and Poonen \cite{OdPo93} observed fractal-like patterns in plots of algebraic integers given by polynomials with coefficients of 0 or 1, using these plots as a jumping off point to prove a number of interesting theorems about the set of such numbers.  Later, Christiensen \cite{Ch06} produced an extremely detailed plot of algebraic integers from a larger coefficient range, coloring the plane according to densities of algebraic integers.  What followed was an even higher resolution image by Sam Derbyshire, inspiring an article in Scientific American \cite{SciAm}.  These images were studied and disseminated further by Baez \cite{Baez11}, who commented on their analytic properties and fractal-like behaviors, and drew connections to Julia sets.

More recent work of Harris, Stange, and Trettel \cite{hst21} serves as the inspiration for the current work and its title.  In their study, they give these images the moniker \emph{algebraic starscapes}.  Rather than color the plane by density of algebraic integers, they plot the integers using dots which are sized inversely proportionally to the \emph{discriminant} of the associated polynomial (which, loosely speaking, serves as a measure of the algebraic complexity of the roots).  This exhibits sharp new visual features for which they offer detailed geometric explanations.  In this way they connect the visual patterns emerging in these \textit{starscapes} to Diophantine approximation and hyperbolic geometry, providing new perspectives to known results and establishing interesting new results along the way.  A similar shift is achieved by StackExchange user DumpsterDoofus \cite{dumpsterdoofus}, who creates plots where (roughly) the points are sized inversely to the size of the coefficients of the associated polynomials, creating stunning visual effects.  Both of these paradigms have the effect of fading the more \emph{algebraically complex} algebraic integers into the background, and the result brings a new level of sharpness and clarity into the picture, highlighting new geometric relationships between algebraic numbers.

A more complete and detailed history of algebraic integer illustrations and art, complete with many references, can be found in \cite[Section 1.4]{hst21}.  Nevertheless, we highlight the paradigm shifts that a change in perspective of the same plot can bring, mathematically and artistically.  While density colorings (\cite{OdPo93},\cite{Ch06},\cite{Baez11}) illuminate analytic and fractal behavior and connections to dynamics, the weighting by discriminant \cite{hst21} brings connections to hyperbolic geometry, and reflects slices of higher dimensional phenomanae, and the weighting by coefficients \cite{dumpsterdoofus} brings forth sharp, enchanting, and mysterious imagery and patterns.  This inspired us to investigate how attaching other sorts of arithmetic invariants to algebraic numbers can highlight other patterns lurking within these \emph{algebraic starscapes}.

Our focus in this paper is to create plots of algebraic integers where the points are sized according to how \emph{symmetric} they are$-$more precisely, according to the size of the \emph{Galois groups} of the associated field extensions.  Gallagher \cite{Gal73} shows that the Galois group of a field extension associated to an irreducible polynomial over the rationals of degree $n$ is almost always the symmetric group $S_n$ (with more precise asymptotics established by Cohen \cite{Coh79},\cite{Coh81}).  Nevertheless, in many cases of interest, the solutions to such a polynomial will satisfy extra algebraic relations, which in turn cuts down the size of the Galois group.  These are the cases that our techniques amplify, showing that these \emph{more rigid} algebraic integers form interesting and beautiful subspaces of the starscapes of all algebraic integers.

The structure of the paper follows.  In Section \ref{GaloisCrash}, we set up the theoretical backbone of our work.  We begin by giving an informal description of the symmetries we are hoping to capture, and then establish the precise language with a brief summary of the underlying Galois theory, introducing the notion of \emph{rigidity}.  We then describe the algorithm for computing rigidity, implementing algorithms of Paulsen \cite{Pau99} to reduce the problem of computing sizes of Galois groups to factorizations over the rational numbers.  This allows us to compute the datasets from which our images are based.  None of the results appearing in Section \ref{GaloisCrash} are new, but as far as the authors could tell, the notion of rigidity as defined here does not seem to appear in the literature.  We therefore include proofs of what are almost certainly known results, but repackaged and clarified in terms of rigidity.

In Section \ref{imageGeneration} we describe how create images from a collection of algebraic integers together with their rigidity, adapting the techniques of DumpsterDoofus \cite{dumpsterdoofus}.  In Section \ref{gallery} we present a selection of the generated images, together with analysis of the structures our methods illuminate.  We compare illustrations of rigid algebraic integers with images sized according to discriminants or coefficients, and interpolations therein.  We also give a complete explanation of the patterns arising in the images of rigid cubics, as well as a partial explanation of the patterns arising in the images of rigid quartics.  In particular, Proposition \ref{RatGeo} shows that rigid quartics concentrate on the so-called \textit{rational geodesics} introduced in \cite{hst21}, which perhaps surprisingly connects the patterns arising in our Galois theoretic analysis of quartics to the patterns arising in quadratics studied in \cite{hst21}.

One can find more images, some generated datasets, and some example code in our public github repository \cite{github}.
\subsection{Acknowledgements}
The core of this project was completed by Xu as part of a senior honors thesis in mathematics at the University of California, Berkeley, under the supervision of Dorfsman-Hopkins.  Xu thanks the Berkeley department of mathematics for their support through the honors program.  The authors thank Katherine Stange for comments on a previous draft, and in particular for suggesting the connection between rational geodesics and the patterns arising in our images of rigid quartics, leading to Proposition \ref{RatGeo}. We also thank the anonymous referees for their detailed and insightful comments.  We could not have generated the images at the center of this manuscript without the help of William Paulsen, who communicated with us on how to use their Mathematica package to efficiently compute sizes of Galois groups of large sets of algebraic integers, as well as StackExchange user DumpsterDoofus who shared code and techniques central to the creation of the images in this paper.  We also must thank Edmund Harris and Steve Trettel with whom Dorfsman-Hopkins shared many inspiring and illuminating conversations while in residence at ICERM.  Indeed, the inspiration for this work can be traced in no small part back to the Fall 2019 semester program at ICERM entitled \emph{Illustrating Mathematics}, which was supported by NSF grant DMS-1439786.  Dorfsman-Hopkins is currently supported by NSF grant DMS-1646385 as part of the Research Training Group in arithmetic geometry at the University of California, Berkeley.
\section{Algebraic Symmetries and Rigid Algebraic Integers}\label{GaloisCrash}
If we plot the roots of a polynomial, we can observe that these solutions obey certain symmetries.  For example, the roots of $x^2+1$ are the imaginary numbers $i$ and $-i$, which satisfy reflection symmetry over the real axis in the complex plane.  This extends to the operation of \emph{complex conjugation}:
\[(z = a+bi) \mapsto (\overline z = a-bi).\]
Importantly, this symmetry is in fact an algebraic operation on the complex plane, in that it respects multiplication and addition:
\[\overline z+\overline w = \overline{z+w}\hspace{15pt}\text{and}\hspace{15pt}\overline z\cdot\overline w = \overline{z\cdot w}.\]
These are the types of symmetries at the heart of our illustrations, the \emph{algebraic symmetries} which respect multiplication and addition.  The general question in Galois theory is then the following: what sorts of algebraic symmetries exist among the roots of a polynomial?  Notice that if we call the roots of a polynomial $r_1,\cdots,r_d$, then an algebraic symmetry of the roots gives a permutation of the set $\{r_1,\cdots,r_d\}$.  A natural question then is whether every permutation of the roots corresponds to an algebraic symmetry.  Gallagher shows \cite{Gal73} that for almost every polynomial, the answer is yes.  Nevertheless, it is not always the case, as the following example illustrates.
\begin{Example}\label{unityExample}
  The roots of $f(x) = x^5-1$ are $\{1,\zeta,\zeta^2,\zeta^3,\zeta^4\}$ where $\zeta=e^{2\pi i/5}$.  Any permutation $\sigma$ of the roots which swaps $\zeta$ and $\zeta^2$ is not algebraic, because it does not respect multiplication:
  \[\sigma(\zeta^2) = \zeta,\hspace{15pt}\text{but}\hspace{15pt}\sigma(\zeta)^2 = (\zeta^2)^2 = \zeta^4.\]
  Therefore not every permutation of the roots of $f$ corresponds to an algebraic symmetry.  The roots of $f$ in this example satisfy an extra algebraic relationship (being powers of each other).  This gives the set of roots some extra rigidity, and therefore limits the amount of algebraic symmetries that can exist among them.  In what follows, we will call such an algebraic integer \emph{rigid}.
\end{Example}
\subsection{Field Extensions and Galois Groups: A Crash Course}
We have so far been rather imprecise about what we mean when we say \emph{algebraic symmetry}.  The set of what we've been calling algebraic symmetries of the roots of a polynomial is more commonly called the \emph{Galois group} of the \emph{splitting field} of a polynomial.  We briefly review the relevant definitions and results, and refer the interested reader to one of the many introductory texts on the subject of Galois theory (for example \cite[Chapters 13 and 14]{DF04}).

At its core, Galois theory is the study of symmetries of field extensions, and much of the numerical data we use in generating our images comes from computing degrees of field extensions.  We begin this section by recording some of the numerical results we will be using throughout.

\begin{Proposition}[\cite{DF04} 13.1/13.2]\label{DegreeFacts}
  Let $K$ be a field and $L$ an extension of $K$.
  \begin{enumerate}
    \item{
    $L$ is a $K$-vector space, whose dimension we call the \textbf{\emph{degree}} of $L$ over $K$ and denote by $[L:K]$.
    }
    \item{
    Degree is multiplicative.  That is, if $F\supseteq L\supseteq K$ are field extensions then:
    \[[F:K] = [F:L][L:K].\]
    }
    \item{
    Let $f(x)\in K[x]$ be an irreducible polynomial of degree $d$.  Then:
    \[L:=K[x]/(f)\]
    is a field extension of degree $d$.  Indeed, $1,x,x^2,\cdots,x^{d-1}$ forms a basis of $L$.}
  \end{enumerate}
\end{Proposition}
Notice that in case 3 above, the polynomial $f$ has a root in $L$ (given by the image $\alpha$ of $x$ in the quotient).  We will therefore write $L = K(\alpha)$, and think of it as the field given by formally adjoining a root of $f$.  We are interested in fields which contain all the roots of a polynomial.
\begin{Definition}
  Let $K$ be a field, and $f(x)\in K[x]$ an irreducible polynomial.
  \begin{enumerate}
    \item{Let $L$ be a field containing $K$.  We say that $f$ \textbf{\emph{splits}} in $L$ if $f$ factors into linear factors as a polynomial in $L[x]$.}

    \item{If $f$ splits in $L$ and there is no strict subfield of $L$ in which $f$ splits, then $L$ is called the \textbf{\emph{splitting field}} of $f$.  These exist and are unique up to isomorphism (\cite[Corollary 13.4.28]{DF04})}.

    \item{Let $L$ be the splitting field of $f$.  The \textbf{\emph{Galois group}} of $L/K$ is the set
    \[\Gal(L/K) = \{\text{isomorphisms }\varphi:L\to L\text{ such that }\varphi(t)=t\text{ for every }t\in K\},\]
    which is a group under composition.  When $K=\mathbb{Q}$, we will call this the Galois group of $f$.}
  \end{enumerate}
\end{Definition}
Notice that a polynomial $f$ splits in $L$ if and only if all of its roots are in $L$, so that the splitting field of $f$ can be thought of as the smallest extension containing all of its roots.  We remark that rigidity will be related to the sizes of Galois groups, and for our purposes, the splitting field will provide us with the numerical means to compute this size.
\begin{Proposition}[\cite{DF04} Proposition 14.1.5]\label{SplittingNumerics}
  Suppose $K$ contains $\bQ$.  If $L$ is the splitting field of an irreducible polynomial in $K[x]$, then:
  \[|\Gal(L/K)| = [L:K].\]
\end{Proposition}
Let $L/K$ be the splitting field of an irreducible polynomial $f$ of degree $d$, and fix $\sigma\in\Gal(L/K)$.  Then, as described above, $\sigma$ gives a permutation of the $d$ roots of $f$, so if we let $S_d$ be the associated permutation group we obtain a homomorphism:
\[\varphi:\Gal(L/K)\to S_d.\]
In fact, $L$ is generated algebraically over $K$ by the roots of $f$ \cite[Section 13.4]{DF04}, so that $\sigma\in\Gal(L/K)$ is completely determined by its values on the roots of $f$.  This means that $\varphi$ must be injective.  Furthermore, $\Gal(L/K)$ acts transitively on the roots of $f-$that is, for any pair of roots of $f$ there must be some $\sigma\in\Gal(L/K)$ taking one to the other$-$so that $\Gal(L/K)$ must have at least $d$ elements.  To summarize, we have shown the following fact:
\begin{Proposition}\label{inequalities}
  Let $f\in K[x]$ be an irreducible polynomial of degree $d$, and let $L/K$ be its splitting field.  Then:
  \[d\le |\Gal(L/K)|\le d!\]
\end{Proposition}
The notion of rigidity has to do with how far the order of the Galois group is from the upper bound.  A result of Gallagher says it's usually not very far.
\begin{Proposition}[\cite{Gal73}]
  Let $N,d\ge 1$.  Let $P_{N,d}$ be the set of polynomials in $\bZ[x]$ whose coefficients have absolute value $\le N$.  Let $S_{N,d}\subseteq P_{N,d}$ be the subset consisting of polynomials whose Galois group maps isomorphically to the full symmetric group $S_d$ via $\varphi$.  Then:
  \[\lim_{N\to\infty}\frac{|S_{N,d}|}{|P_{N,d}|} = 1.\]
\end{Proposition}
This says that the situation in Example \ref{unityExample} exhibits a rather rare phenomenon, which is what we call \emph{rigidity}.  Before we can define this precisely, we must describe how to assign a Galois group to an algebraic integer.  We will need the following fact.
\begin{Proposition}[\cite{DF04} Section 13.2]\label{MinimalPoly}
  Let $\alpha$ be the solution to some polynomial over $\bQ$.  Then there is a unique monic irreducible polynomial $m_\alpha\in\bQ[x]$ which has $\alpha$ as a root.
\end{Proposition}
This allows us to define the objects at the center of our study.
\begin{Definition}
  If $\alpha$ is the solution to some polynomial over $\bQ$, we call $\alpha$ an \textbf{\emph{algebraic number}}, and we call the polynomial $m_\alpha$ from Proposition \ref{MinimalPoly} the \textbf{\emph{minimal polynomial}} of $\alpha$.  The \textbf{\emph{degree}} of $\alpha$ is the degree of $m_\alpha$, and the \textbf{\emph{Galois group}} of $\alpha$ is the the Galois group of $m_\alpha$ (that is, $\Gal(L/\bQ)$ where $L$ is the splitting field of $m_\alpha$).  We call $\alpha$ an \textbf{\emph{algebraic integer}} if the $m_\alpha$ has coefficients in $\bZ$.
\end{Definition}
\begin{Example}\label{classicalIntegers}
  The algebraic integers of degree one are precisely the integers.  Indeed, the minimal polynomial of such an $\alpha$ would have to be $x-\alpha\in\bZ[x]$, making $\alpha$ an integer.  The Galois group of $\alpha$ is the trivial group.
\end{Example}
\begin{Example}\label{quadraticIntegers}
  Let $n$ be an integer which is not a perfect square.  Then $\sqrt n$ is an algebraic integer of degree $2$, whose minimal polynomial is $x^2-n$.  The Galois group of $\sqrt n$ is isomorphic to the cyclic group $\bZ_2$, whose nontrivial element corresponds to \emph{conjugation}: $\sqrt n\mapsto -\sqrt n$.
\end{Example}
\begin{Example}\label{unityExample2}
  The element $\zeta$ from Example \ref{unityExample} is an algebraic integer, but its minimal polynomial is \emph{not} $f(x) = x^5-1$, since it is not irreducible.  Instead, the minimal polynomial is the cyclotomic polynomial $m_\zeta(x) = x^4+x^3+x^2+x+1$, which we obtain after factoring out $x-1$ from $f$.  Therefore $\zeta$ is an algebraic integer of degree $4$.  We next compute its Galois group $G$.  Notice that $\sigma\in G$ is determined by where it sends $\zeta$, since the other roots of $m_\zeta$ are all powers of $\zeta$ and $\sigma$ must respect multiplication.  That is, $\sigma$ is determined by $\sigma(\zeta) = \zeta^i$ for some $i=1,2,3,4$, and so we conclude that $G\cong \bZ_4$.  Notice that this is much smaller than $S_4$ which has order $4! = 24$.
\end{Example}
\subsection{Rigidity: Definition and Basic Properties}
We can now define the main invariant we are illustrating in this project.
\begin{Definition}\label{rigidityDef}
  Let $\alpha$ be an algebraic integer of degree $d$ with Galois group $G$.  The \textbf{\emph{rigidity}} of $\alpha$ is:
  \[\operatorname{rig}(\alpha):=\frac{d!-|G|}{d!}\]
  If $\rig(\alpha)>0$ then we call $\alpha$ \textbf{\emph{rigid}}.
\end{Definition}
Notice that an immediate consequence of Proposition \ref{inequalities} is that
\[0\le\rig(\alpha)\le1-\frac{1}{(d-1)!}.\]
This means that $\alpha$ never has rigidity equal to 1.  This is reasonable in a qualitative sense, since something that is fully rigid should admit no nontrivial automorphisms, but algebraic integers (of degree $\ge2$) will always admit some symmetries.  If $\alpha$ is an algebraic integer of the maximal allowable rigidity, then we call $\alpha$ \emph{maximally rigid}.
\begin{Example}\label{unityExample3}
  Consider $\zeta$ from Examples \ref{unityExample} and \ref{unityExample2}.  Our computations show that:
  \[\rig(\zeta) = \frac{4!-4}{4!} = \frac{5}{6},\]
  so that $\zeta$ is maximally rigid.
\end{Example}
\begin{Example}\label{fourthRootOf2}
  We now compute the rigidity of $\sqrt[4]{2}$, which demonstrates how one might compute rigidity in general.  The idea is to compute the degree of the splitting field of the minimal polynomial $f(x) = x^4-2$, and then leverage Proposition \ref{SplittingNumerics}.  We first consider $K = \bQ(\sqrt[4]{2}) = \bQ[x]/(f(x)),$ which is a field extension of degree 4 by Proposition \ref{DegreeFacts}(3).  Then in $K$, $f(x)$ factors as
  \[x^4-2 = (x-\sqrt[4]{2})(x+\sqrt[4]{2})(x^2+\sqrt{2}).\]
  Letting $g(x) = x^2+\sqrt{2}\in K[x]$, we define $L = K[x]/(g(x)) = K(i\sqrt[4]{2})$, which is an extension of degree 2.  Then in $L$, we observe that $f$ splits completely as:
  \[f(x) = (x-\sqrt[4]{2})(x+\sqrt[4]{2})(x-i\sqrt[4]{2})(x+i\sqrt[4]{2}).\]
  Therefore $L$ is the splitting field of $f(x)$.  Then by applying Propositions \ref{SplittingNumerics} and \ref{DegreeFacts}(2), we compute the size of the Galois group of $\sqrt[4]{2}$ as:
  \[|\Gal(L/\bQ)| = [L:\bQ] = [L:K][K:\bQ] = 2\times 4 = 8.\]
  We can therefore compute that $\rig(\sqrt[4]{2}) = \frac{2}{3}$.  In particular, we see that it is rigid, but not maximally rigid.  Identical considerations hold for the remaining roots of $f$.
\end{Example}
Before moving forward, we record some easily deduced facts about rigidity for algebraic integers of degree less than or equal to 3.
\begin{Proposition}\label{le2}
  Algebraic integers of degree $\le2$ are never rigid.
\end{Proposition}
\begin{Proof}
  This is an immediate consequence of Proposition \ref{inequalities}, since for $d=1,2$, we have $d!=d$, so that the order of the Galois group would be equal to both, and therefore the rigidity is easily computed to be 0.
\end{Proof}
This says that for algebraic integers of degree 1 or 2, the notion of rigidity adds nothing new to the picture, so we should not expect the images generated using our methods to provide new insights for these small degrees.  Rigid algebraic integers of degree 3 do exist, but the following fact tells us that they are all concentrated on the real axis of the complex plane.
\begin{Proposition}\label{cubics}
  If an algbraic integer of degree 3 is rigid, it has rigidity equal to $\frac{1}{2}$ so that it is maximally rigid.  Furthermore, it is a real number.
\end{Proposition}
\begin{Proof}
  Let $f\in\bQ[x]$ be an irreducible polynomial of degree 3, let $L$ be its splitting field, and let $G = \Gal(L/\bQ)$.  By Proposition \ref{inequalities}, we know that $3\le|G|\le 6$, and further, by Lagrange's theorem we know $|G|$ divides $|S_3| = 6$, so it can only take the values $3$ or $6$.  Therefore if $f$ is the minimal polynomial of a rigid algebraic integer, we must have $|G|=3$.  The first statement then follows from the definitions.  Now let $K = \bQ(\alpha) = \bQ[x]/(f(x))$ be the subfield of $L$ generated by a single root $\alpha$ of $f$.  Then by Proposition \ref{DegreeFacts}(2) we observe that
  \[3 = [L:\bQ] = [L:K][K:\bQ] = [L:K]\cdot3,\]
  so that $[L:K]=1$ and $L=K$.  In particular, the splitting field of $f$ is generated by adjoining any single (possibly complex) root of $f$ to $\bQ$.  But by the intermediate value theorem, a cubic polynomial has at least one real root, say $\tau\in L\cap\bR$, so that $L = \bQ(\tau)$ is algebraically generated by $\bQ$ and $\tau$.  However, any algebraic combination of real numbers must be real, so $L\subseteq\bR$, and in particular every root of $f$ must have been real to begin with.
\end{Proof}
The qualitative consequence of Proposition \ref{cubics} is that any modifications we make to images of degree 3 algebraic integers using the notion of rigidity will be concentrated along the real axis, meaning they probably won't make very interesting patterns in $\bC$.
\begin{Remark}\label{RigidReality}
  Notice that the proof of Proposition \ref{cubics} showed more generally that if the Galois group of a polynomial $f$ of degree $d$ has $d$ elements, and $f$ has at least one real root, then all of the roots of $f$ are real.  Since all polynomials of odd degree have a real root, this says that any maximally rigid algebraic integer of odd degree is a real number.  Qualitatively, this suggests that we should expect concentrations of rigid numbers along the real axis, especially in the odd degree cases.
\end{Remark}
\subsection{Computing Rigidity}
We now describe the algorithm we implement to compute the rigidity of an algebraic integer, essentially following the steps of Example \ref{fourthRootOf2}.  It boils down to following the standard construction of the splitting field of a polynomial, while keeping track of some extra data along the way.
\begin{Algorithm}\label{bigAlg}
  Given a field $K$ and a polynomial $f\in K[x]$, the following computes the splitting field of $f$ as well as its degree.
  \begin{enumerate}
    \item{Set $i = 0$, $K_0 = K$, and $n = 1$.}
    \item{Factor $f(x) = f_1(x)f_2(x)\cdots f_r(x)\in K_i[x]$ into irreducibles.}
    \item{If all the $f_j$ are linear, return $n$ and $K_i$.  Else reorder the factors so that $f_1$ has degree $d>1$.}
    \item{Set $K_{i+1} = K_i[x]/f_1(x)$, $f = f/f_1$ and $n = nd$.  Return to step 2.}
  \end{enumerate}
\end{Algorithm}
\begin{nProof}[Proof of Correctness]
  The fact that the algorithm terminates is clear as the degree of $f$ decreases at each step.  The fact that it computes the splitting field is \cite[Theorem 13.4.25]{DF04}.  The fact that it correctly computes the degree of the splitting field follows by considering the inclusions
  \[K = K_0\subseteq K_1\subseteq\cdots\subseteq K_t,\]
  where $K_t$ is the splitting field, and applying Proposition \ref{DegreeFacts}(2) inductively.
\end{nProof}
Once we have calculated the degree of the splitting field, we have (by Proposition \ref{SplittingNumerics}) computed the size of the Galois group, and therefore can easily determine the rigidity of all the roots of $f$ by just plugging this value into the formula of Definition \ref{rigidityDef}.

To implement step 2 of Algorithm \ref{bigAlg}, we must factor a polynomial over an extension of $\bQ$.  Our implementation adapts a Mathematica package of Paulsen \cite{Pau99} called \verb|galois.m|, which uses the theory of field norms to reduce the factorization to one over the rational numbers, where there are well known (and fast) algorithms.
\begin{Definition}\label{normDef}
  Let $f(x)\in\bQ[x]$ be a polynomial with (complex) roots $\alpha = \alpha_1,\cdots,\alpha_n$, and let $g = g(\alpha,x)\in \bQ(\alpha)[x]$.  We define the \textbf{\emph{norm}} of $g$ to be:
  \[N(g(\alpha,x)) = g(\alpha_1,x)g(\alpha_2,x)\cdots g(\alpha_n,x).\]
\end{Definition}
\begin{Lemma}
  In the setup of Definition \ref{normDef}, we have $N(g(\alpha,x))\in\bQ[x]$.
\end{Lemma}
\begin{Proof}
  Let $L$ be the splitting field of $f$, so that in particular each $\alpha_i\in L$.  Then $N(g(\alpha,x))\in L[x]$.  Since any $\sigma\in\Gal(L/\bQ)$, permutes the $\alpha_i$ and leaves elements of $\bQ$ fixed, we see $\sigma$ just permutes the factors of $N(g(\alpha,x))$ so that $\sigma(N(g(\alpha,x)))=N(g(\alpha,x))$.  Therefore $N(g(\alpha,x))$ has coefficients in $\bQ$ by the Galois correspondence \cite[Theorem 14.2.14]{DF04}.
\end{Proof}
It turns out that the factorization of $g$ in $\bQ(\alpha)[x]$ is closely related to the factorization of $N(g)$ in $\bQ[x]$.  Indeed, since the norm is multiplicative, a factorization of $g$ would give a factorization of $N(g)$.  Unfortunately, the reverse direction isn't true in general, but if we modify $g$ slightly before taking the norm we do obtain a converse, as the following proposition of Paulsen \cite[Propostition 2.2]{Pau99} shows.  We reproduce the proof and fill in some details for the convenience of the reader.
\begin{Proposition}\label{Norming}
  In the setup of Definition \ref{normDef}, define $h(\alpha,x,y):=g(\alpha,xy-\alpha y)\in\bQ(\alpha)[x,y]$ where $y$ is a new variable.  Then $g(\alpha,x)$ is irreducible in $\bQ(\alpha)[x]$ if and only if $N(h(\alpha,x,y))$ is irreducible in $\bQ[x,y]$.
\end{Proposition}
\begin{Proof}
  The backward direction follows immediately from the multiplicativity of the norm.  To prove the forward direction we may assume without loss of generality that $g(\alpha,x)$ is monic and irreducible of degree $r$.  Consider the factorization:
  \[N(h(\alpha,x,y)) = g(\alpha_1,xy-\alpha_1 y)\cdots g(\alpha_n,xy-\alpha_n y)\in L[x,y],\]
  where $L$ is a splitting field for $f$.  Denote by $g_i:=g(\alpha_i,xy-\alpha_i y)\in\bQ(\alpha_i)[x,y]$.  We will first prove that since $g$ is irreducible, each $g_i$ is.  Indeed, by symmetry it suffices to prove $g_1$ is irreducible.  But given a factorization of $g_1$, setting $y=1$ gives a factorization of $g(\alpha,x-\alpha)$, which by a change of coordinates factors $g$.

  Now suppose $N(h)$ factors as $pq\in\bQ[x,y]$.  Then in $\bQ(\alpha_i)[x,y]$ we have $pq = g_i l_i$, so that by the irreducibility of $g_i$, we know $g_i$ divides $p$ or $g_i$ divides $q$.  Doing this for each $i$, and recalling that $g$ is monic of degree $r$, shows that $p$ has the following form:
  \begin{eqnarray*}
    p(x,y) &=& g(\alpha_{i_1},xy-\alpha_{i_1} y)\cdots g(\alpha_{i_t},xy-\alpha_{i_t}y)\\
    &=&y^{tr}(x-\alpha_{i_1})^r(x-\alpha_{i_2})^r\cdots(x-\alpha_{i_t})^r + \text{lower order terms in }y.
  \end{eqnarray*}
  In particular, we see that $(x-\alpha_{i_1})^r(x-\alpha_{i_2})^r\cdots(x-\alpha_{i_t})^r\in\bQ[x]$, which implies (since $\bQ$ is a perfect field) that $(x-\alpha_{i_1})(x-\alpha_{i_2})\cdots(x-\alpha_{i_t})\in\bQ[x]$ as well.  But this clearly divides $f(x)$ (seeing as the $\alpha_i$ are precisely the roots of $f$), so we must have that $t=0$ or $t=n$ as $f$ is irreducible in $\bQ[x]$.  In particular, $p=N(h)$ or $p=1$, and so $N(h)$ was irreducible to begin with.
\end{Proof}
Proposition \ref{Norming} suggests that in order to factor $g(x)$ over $\bQ(\alpha)[x]$, we may instead introduce a variable and take norms to factor $N(g(xy-\alpha y))\in\bQ[x,y]$, reducing it to a question over the rational numbers.  What remains is to recover the factorization of $g$ from the factorization of $N(g(xy-\alpha y))$.  This boils down to the following question: suppose we are given a polynomial $p(x,y)\in\bQ[x,y]$ which is known to be the norm $N(q(xy-\alpha y))$ for some polynomial $q(x)\in \bQ(\alpha)[x]$, can we recover $q(x)$?  In fact, Paulsen's Mathematica package \verb|galois.m| has a built in function \verb|unNorm| that does exactly this.  The general idea is the following: first notice the degree of $q(x)$ is the product of the degree of $p(x)$ and the degree of the extension.  Then take a generic polynomial $F(x)\in \bQ(\alpha)[x]$ of that degree, and compute $N(F(xy-\alpha y))$.  This way one obtains a system of equations by considering the various coefficients of $y$, which can be efficiently solved for the coefficients of $F$.  Documentation for \verb|galois.m| can be found in \cite{Pau99}.  A version of the package has been patched by the second author to run with the current version of Mathematica, and can be found in our repository \cite{github}.
\section{Generating Images}\label{imageGeneration}
 In this section we will explicitly describe the process of turning the information gathered in Section \ref{GaloisCrash} into our desired \emph{starscapes}.  Since algebraic integers in the complex plane are dense, one could color the plane according to the density of the algebraic integers.  As discussed in the introduction, this was the technique of Derbyshire and others. Alternatively, one could plot the algebraic integers as points, and adjust emphasis by sizing these points according to various criteria.  For example, \cite{dumpsterdoofus} sizes points inversely proportionally to the size of the coefficients of the minimal polynomial, which tends to shrink points with more complicated equations.  Qualitatively, this highlights the more \emph{algebraically simple} integers.  Similarly, \cite{hst21} sizes points inversely to the size of the discriminant, which again highlights more algebraically simple integers (for a slightly different definition of simple).

Our method is closer to the second paradigm, but in contrast, we enlarge points which are more rigid, to highlight the geometric relationships that rigid integers hold.  We also generate images using discriminant or coefficient weightings as in \cite{hst21} and \cite{dumpsterdoofus}, as well as interpolations between those and the rigidity weightings, and compare the results.

Once we have gathered a collection of algebraic integers together with rigidity and other invariants of interest we choose a brightness function.  This function will assign a brightness to each algebraic integer in terms of the invariants of interest, so that in the final plots integers with higher brightness values will be emphasized, appearing brighter and larger.  Our overall process boils down to the following 3 steps.
\begin{enumerate}
\item{
Generate a dataset of algebraic integers together with the invariants of interest.  In particular, the elements will consist of tuples $(z,r,|v|,D)$ where $z$ is an algebraic integer with rigidity $r$ and discriminant $D$, and $|v|$ is the magnitude of the coefficient vector $(1,a_{d-1},\cdots,a_0)$ where $x^d + a_{d-1}x^{d-1} + \cdots + a_0$ is the minimal polynomial of $z$.
}
\item{
Choose a brightness function $\beta: = \beta(z,r,|v|,D)$ which assigns a brightness to each algebraic integer in terms of rigidity, discriminant, and the magnitude of the coefficient vector.
}
\item{
Create a plot of a region in the complex plane satisfying the following goal: for each algebraic integer $z = a+bi$ in our list, draw a dot at the coordinates $(a,b)$, whose size and brightness is controlled by the brightness value $\beta$ associated to $z$.
}
\end{enumerate}
Let us go into slightly more detail about how we achieve each step.

For the first step we enumerate a list of monic \emph{irreducible} polynomials $x^d+a_{d-1}x^{d-1}+\cdots+a_0$ of a given degree, where the $a_i$ are integers within some coefficient bound $[-c,c]$ for some integer $c$.  At this point we also record the magnitude of the coefficient vector $v = (1,a_{n-1},\cdots,a_0)$, as well as the discriminant $D$ of the polynomial, which are both computed algebraically in terms of the coefficients.  We then compute the roots of each polynomial (which will be degree $d$ algebraic integers), together with their rigidity $r$, using the methods of Section \ref{GaloisCrash}.  Our data sets are in a public repository \cite{github}.
\begin{table}[h!] 
  \centering
  \begin{tabular}{|c:c:c|}
    \hline
    $\beta_D$ & ${1}/{\sqrt{|D|}}$ &  Small discriminant emphasis\\
    \hline
    $\beta_v$ & ${1}/{|v|}$ & Small vector magnitude emphasis\\
    \hline
    $\beta_r$ & $(1+2\rig)^4$ & Rigidity emphasis\\
    \hline
    $\beta_r^+$ & $(1+3\rig)^4$ & Rigidity: extreme emphasis\\
    \hline
  \end{tabular}
  \caption{Brightness parameters}
  \label{brightnessParameters}
\end{table}

For the second step, we explore various brightness functions $\beta$, listed in Table \ref{brightnessParameters}.  We point out that the brightnesses defined in the table are only \emph{proportional} to the formulas given, and when generating images in practice they are scaled by a constant.  We also generate images interpolating between these weightings, using the following brightness formulas:
\begin{eqnarray}\label{compoundBrightness}
  \beta_{Dr} := \beta_D\times\beta_r&\hspace{30pt}&\beta_{vr} := \beta_v\times\beta_r\\
  \beta_{Dr}^+ := \beta_D\times\beta_r^+&\hspace{30pt}&\beta_{vr}^+ := \beta_v\times\beta_r^+.\nonumber
\end{eqnarray}
For example, $\beta_{Dr}^+$ fades integers of larger discriminant into the background, while also enlarging rigid points, with the approximate outcome of emphasizing the more \emph{algebraically simple} rigid points.

For the third step we adopt the techniques of StackExchange user DumpsterDoofus \cite{dumpsterdoofus}.   We first translate the data from steps 1 and 2 into a \emph{sparse matrix}, whose entries should represent the pixels on the screen.  To do this, we start with a large matrix, all of whose entries are 0.  For each tuple $(a+bi,r,|v|,D)$, we scale and round $a$ and $b$ to integers $A$ and $B$.  We then input the brightness $\beta := \beta(a+bi,r,|v|,D)$ into the $(A,B)$ entry of the matrix.   We generate an image from this matrix by coloring the $(i,j)$-pixel of the image according to the brightness $\beta$ of the $(i,j)$-entry of the matrix.  We then apply a Gaussian blur so that brightness translates into the size of a point, and fades rather than appearing as a solid dot.  We refer the interested reader to the StackExchange posts of DumpsterDoofus \cite{dumpsterdoofus} for more detailed descriptions of their strategy, together with detailed example code (which we borrow from heavily).
\begin{Remark}\label{cubeConundrum}
When illustrating mathematics it is important to ask whether the patterns that arise are products of the underlying mathematics, or artifacts the tools used to make the illustration.  In our case, we are plotting integers coming from polynomials whose coefficients lie between $-c$ and $c$. This corresponds to considering a cube in the space of polynomials.  One should ask whether the patterns that appear in the images that follow might reflect the geometry of this cube rather than properties of rigidity. To address this we considered polynomials appearing in a sphere rather than a cube, only plotting integers such that the magnitude of the coefficient vector $|v|$ was bounded by $c$.  A second approach interpolates between the sphere and the cube, and has already been suggested above.  In particular, by using brightness parameters $\beta_{vr}$ or $\beta_{vr}^+$, we shrink points coming from integers whose minimal polynomial is further from the origin in the space of polynomials.   This de-emphasizes the corners of the cube, limiting the effect of its geometry on our overall picture. The patterns we discuss in the following section appear in each case, giving us reasonable certainty that they come from the geometry of rigid numbers rather than the geometry of our underlying space of polynomials.  For the images included in this manuscript we always use the cubical polynomial space, although we do have images using $\beta_{vr}$ in Figures \ref{Frontpage} and \ref{quintics!}.  Find a more complete collection of generated images at \cite{github}.
\end{Remark}

\section{Gallery, Analysis, and Conclusions}\label{gallery}
We begin our discussion of the generated images by considering the four illustrations that appear on the front page of this report (Figure \ref{Frontpage}).  These are generated by fixing a brightness function $\beta$, producing an image following the above steps for irreducible polynomials of degrees 1 through 5, with coefficients in the range $[-4,4]$.  Each degree was colored as in Table \ref{ColoringTable}.

\begin{table}[h!]
  \centering
  \begin{tabular}{|c:c:c|}
    \hline
    Linear & Yellow &\colorbox{black}{$\color[rgb]{.96,.71,.1}\bullet$}\\
    \hline
    Quadratic & Pink&\colorbox{black}{$\color[rgb]{1, .25, .24}\bullet$}\\
    \hline
    Cubic & Purple&\colorbox{black}{$\color[rgb]{.35,.05,.75}\bullet$}\\
    \hline
    Quartic & Blue&\colorbox{black}{$\color[rgb]{.28, .53, .78}\bullet$}\\
    \hline
    Quintic & White&\colorbox{black}{$\color[rgb]{1,1,1}\bullet$}\\
    \hline
  \end{tabular}
  \caption{Colorings in Figure \ref{Frontpage}}
  \label{ColoringTable}
\end{table}

The brightness functions used were $\beta_{D}$ and $\beta_{vr}$ in the top row, and $\beta_r$ and $\beta_{Dr}$ in the bottom row (where we refer to the brightness values in Table \ref{brightnessParameters} and Equation \ref{compoundBrightness}).  One sees that different patterns seem to arise with the different perspectives, producing what we think are intriguing visual effects.  That being said, one should remain cautious when drawing too many conclusions about relationships of integers of different degrees in these images, since the invariants in question behave very differently across different degrees.  For example, discriminants are much larger for integers of higher degrees, leading to a prominence of low degree points in the first and fourth images.  On the other hand, by Propositions \ref{le2} and \ref{cubics}, there are no rigid degree 1 or 2 points, and very few rigid degree 3 points, leading to them being almost invisible in the third image.  For this reason, in what follows we will only consider images consisting of algebraic integers of a fixed degree.

\begin{figure}[h!]
  \centering
  \includegraphics[scale=.22]{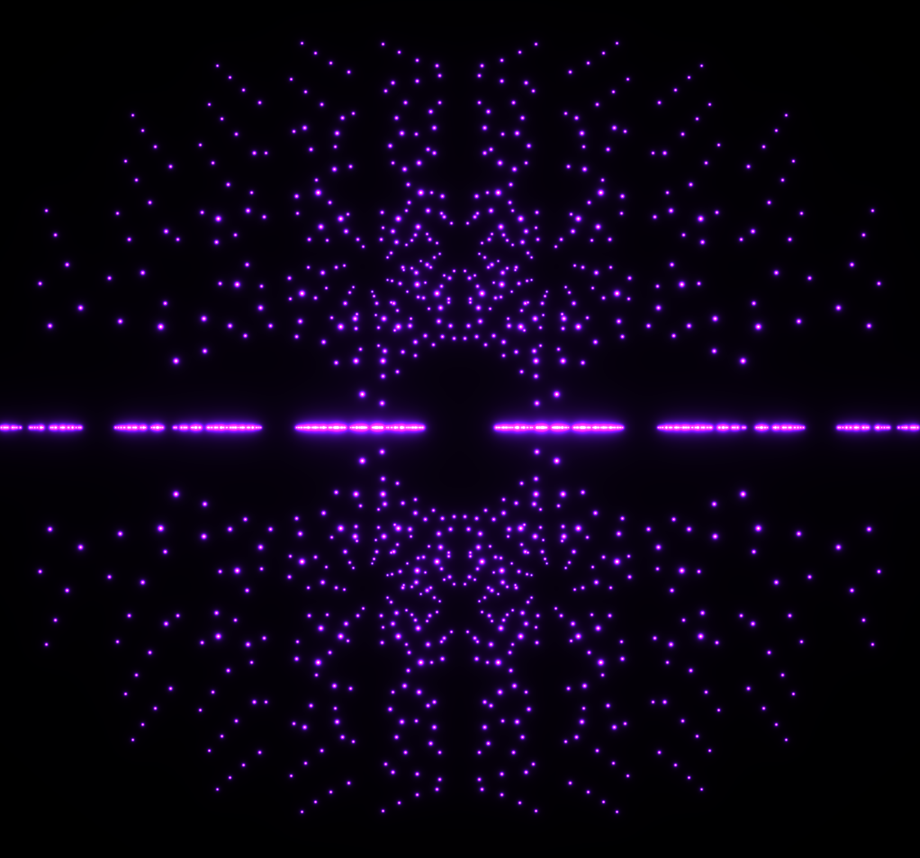}\hspace{20pt}\includegraphics[scale=.22]{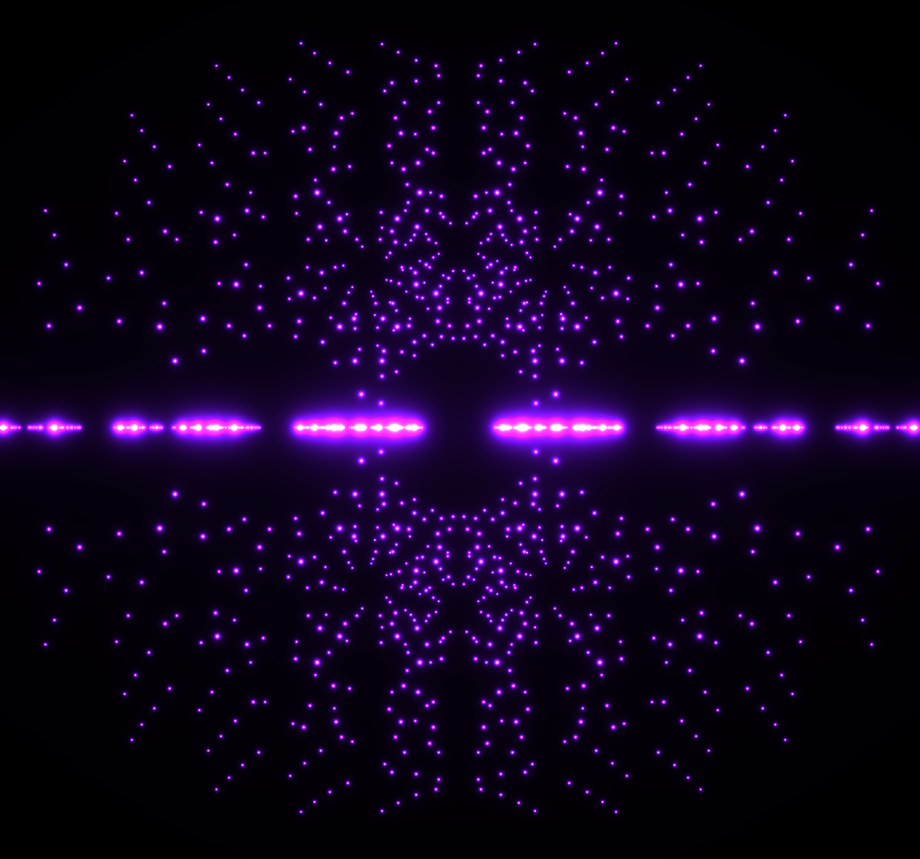}
  \caption{Algebraic integers of degree 3.  \textit{Left}: Sized by discriminant only, using $\beta_D$ as the brightness parameter.  \textit{Right} Sized by disciminant and rigidity, using $\beta_{Dr}$ as the brightness parameter.}
  \label{Cubics!}
\end{figure}

In light of Proposition \ref{le2}, we skip the cases of linear and quadratic integers, since they are never rigid.  We continue our analysis considering case of cubics in Figure \ref{Cubics!}.  Proposition \ref{cubics} tells us exactly what we should expect, that the only rigid cubics appear along the real axis, with rigidity equal to $\frac{1}{2}$.  We use brightness parameters $\beta_D$ and $\beta_{Dr}$ respectively, so that the second is a copy of the first but with rigid points enlarged.  And indeed, it should come as no surprise that the two images are so similar.  As Proposition \ref{cubics} predicts, the second image matches the first except along the real axis, which shines brighter due to the isolated concentration of rigid cubics.

\begin{figure}[h!]
  \centering
  \includegraphics[scale=.22]{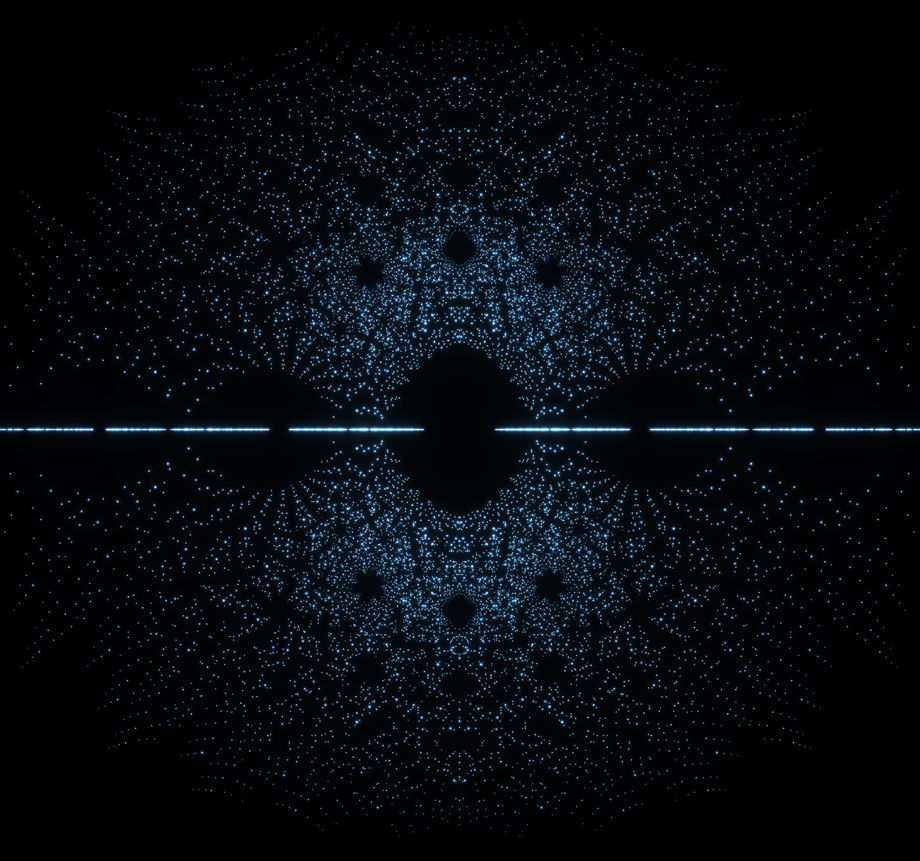}\hspace{20pt}\includegraphics[scale=.22]{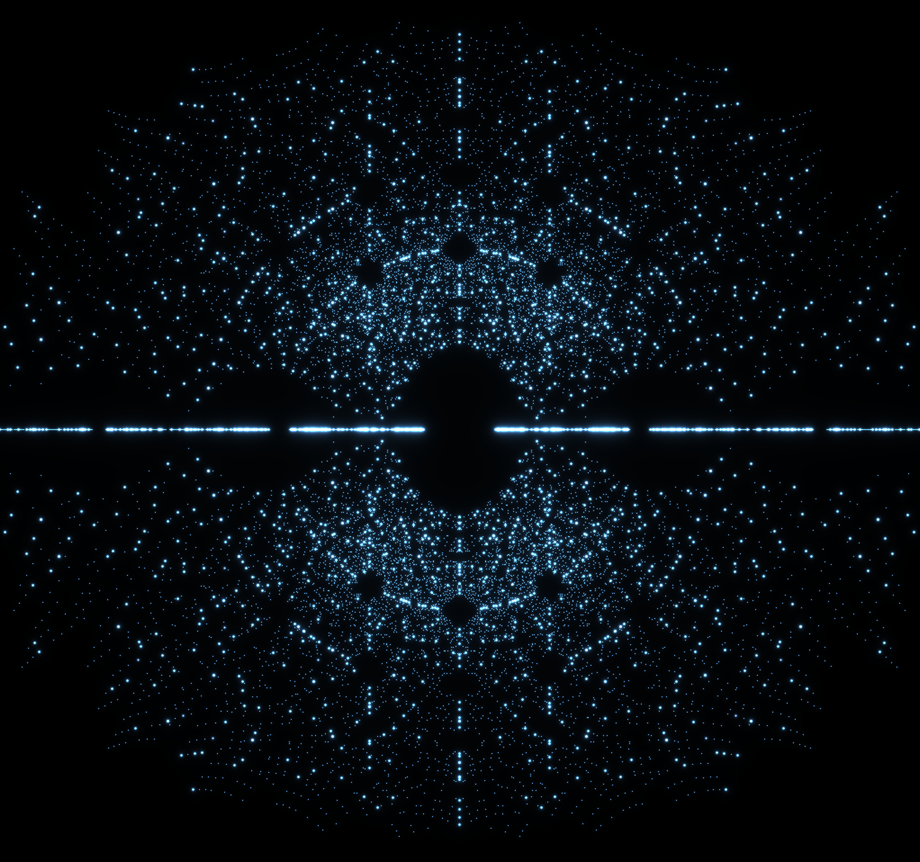}
  \caption{Algebraic integers of degree 4.  \emph{Left}: Sized by discriminant only, using $\beta_D$ as the brightness parameter.  \emph{Right}: Sized by rigidity only, using $\beta_r^+$ as the brightness parameter.}
  \label{quartics!}
\end{figure}

We next consider quartics, depicted in Figure \ref{quartics!}.  Notice that by highlighting rigidity, new patterns emerge in the righthand picture that do not appear in the lefthand one.  Specifically, one can clearly perceive 2 circles and 3 vertical lines.  In fact, these appear to be a partial shadow of the so-called \textit{rational geodesics} introduced in \cite{hst21}.  In particular, in \cite[Figure 1(a)]{hst21}, quadratic points are plotted, sized inversely proportionally to the discriminant.  What stands out are intricate patterns of vertical lines and circles.  Harris, Stange, and Trettel then observe that these lines and circles, which they call \textit{rational geodesics}, correspond precisely to complex numbers $\alpha\in\bC$ satisfying the property that the set $\{1,\alpha\overline \alpha,\alpha+\overline \alpha\}$ is $\bQ$-linearly dependent \cite[Observation 4.12]{hst21}.  If $\alpha$ is a quadratic integer, then a dimension count shows that this dependance must be satisfied, explaining the concentration of quadratic integers on these rational geodesics.

It is perhaps surprising that a pattern arising in plots quadratic integers sized by discriminant arises again in quartic integers sized by the Galois theoretic invariant of rigidity.  In particular, the Galois theory of quadratics plays no role, so it would be notable if applying a Galois theoretic lens to quartics recovers patterns arising in quadratics. Figure \ref{quartics!} suggests it may be true, and the following proposition adds to the evidence, showing that quartics lying on rational geodesics are either maximally rigid, or else quite nearly so, as well as a partial converse.
\begin{Proposition}\label{RatGeo}
Let $\alpha\in\bC$ be a quartic algebraic integer.  If $\alpha$ lies on a rational geodesic, then $\alpha$ is rigid, with rigidity equal to either $\frac{2}{3}$ or $\frac{5}{6}$ (in the latter case $\alpha$ is maximally rigid).  Conversely, if $\alpha$ is maximally rigid with nonzero imaginary part, then it lies on a rational geodesic.
\end{Proposition} 
\begin{Proof}
Let $f$ be the minimal polynomial of $\alpha$, and $L$ the splitting field of $f$.  Notice that the complex conjugate $\overline\alpha$ is a root of $f$ as well, since $f(\overline\alpha) = \overline{f(\alpha)} = \overline0 = 0$.  In particular, $\overline\alpha\in L$.  For the forward direction, we must show that $[L:\bQ]$ is either $4$ or $8$.   We have the following diagram of fields, where solid lines denote field extensions labeled by degree of the extension.
\[
\begin{tikzcd}
&L\ar[dl,dash,swap,"d"]\ar[ddr,dash,"d'"]&\\
\bQ(\alpha)\ar[ddr,dash,swap,"4"]\ar[drr,dash,dashed] &&\\
&& \bQ(\alpha+\overline\alpha)\ar[dl,dash]\\
&\bQ&
\end{tikzcd}
\]
By Proposition \ref{DegreeFacts}(2) we know that $[L:\bQ] =  [L:\bQ(\alpha)][\bQ(\alpha):\bQ] =4d$.  We will deduce the result from the following two observations.  First, the diagonal dotted line can actually be filled in, that is, $\bQ(\alpha+\overline\alpha)\subseteq\bQ(\alpha)$.  Second, $d'=1$ or $2$.  The first observation implies that $d\le d'$, so that it is also 1 or 2, which is exactly what we want.

We now fill in the dotted line.  This will follow if we show that $\overline\alpha\in\bQ(\alpha)$, because then $\alpha+\overline\alpha$ is contained in $\bQ(\alpha)$ and therefore so is the field it generates.  Using that $\alpha$ lies on a rational geodesic, we know by \cite[Observation 4.12(a)]{hst21} that $\{1,\alpha\overline\alpha,\alpha+\overline\alpha\}$ are linearly dependent over $\bQ$.  We then write $c_1 + c_2\cdot\alpha\overline\alpha + c_3(\alpha+\overline\alpha) = 0$ for $c_i\in\bQ$ not all zero, and solve for
\[\overline\alpha = \frac{c_1-c_3\alpha}{c_2\alpha + c_3}\in\bQ(\alpha).\]
It remains to show that $d'$ is either 1 or 2.  Again by the $\bQ$-linear dependence of $\{1,\alpha\overline\alpha,\alpha+\overline\alpha\}$, we see that $\alpha\overline\alpha\in\bQ(\alpha+\overline\alpha)$, so that
\[(x-\alpha)(x-\overline\alpha) = x^2-(\alpha+\overline\alpha)x+\alpha\overline\alpha\in\bQ(\alpha+\overline\alpha)[x].\]
In particular, the minimal polynomial for $\alpha$ over $\bQ(\alpha+\overline\alpha)$ has degree less than or equal to 2, so that the splitting field does as well.   (In fact, since $\alpha+\overline\alpha$ is a real number, the degree is 1 precisely when $\alpha$ is real.)

To prove the converse, we observe that if $\alpha$ is maximally rigid, then $L = \bQ(\alpha)$ (this is essentially the definition of being maximally rigid).  Furthermore, since $\alpha$ has nonzero imaginary part and $\bQ(\alpha+\overline\alpha)$ is totally real, $\bQ(\alpha)\supsetneq\bQ(\alpha+\overline\alpha)\supsetneq\bQ$.  Since the degree of the total extension is 4, we deduce that each intermediate extension must have degree 2.  Because $\bQ(\alpha+\overline\alpha)$ is totally real, it is contained in the fixed field of complex conjugation, but the fixed field of complex conjugation also has degree 2 in $L$, therefore they must be the same.  In particular, since $\alpha\overline \alpha$ is fixed by complex conjugation, it is contained in $\bQ(\alpha+\overline\alpha)$.  But now we see that $1,\alpha+\overline\alpha$, and $\alpha\overline\alpha$ are all contained in this field, which has dimension 2 over the rationals.  Therefore they must be dependent, completing the proof.
\end{Proof}
Proposition \ref{RatGeo} explains the prominence of rational geodesics in the images of rigid quartics, but unlike the case of cubics it is not an exhaustive description of the patterns arising in Figure \ref{quartics!}.  It does almost completely describe the concentration of maximally rigid quartics, but it says nothing about quartics with rigidity equal to $\frac{1}{2}$, other than that they do not lie on rational geodesics.  It also leaves open the possibility quartics of rigidity $\frac{2}{3}$ lying elsewhere as well.  In particular, it begins to paint a picture, but leaves open the possibility of further understanding.

\begin{figure}[h!]
  \centering
  \includegraphics[scale=.22]{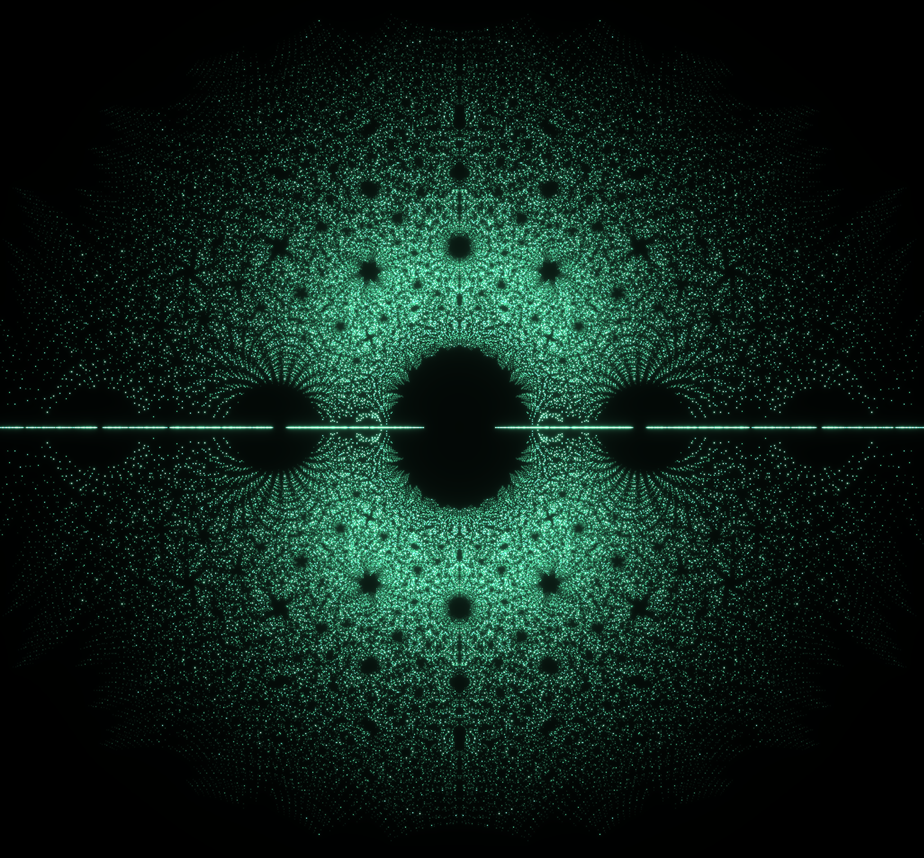}\hspace{20pt}\includegraphics[scale=.22]{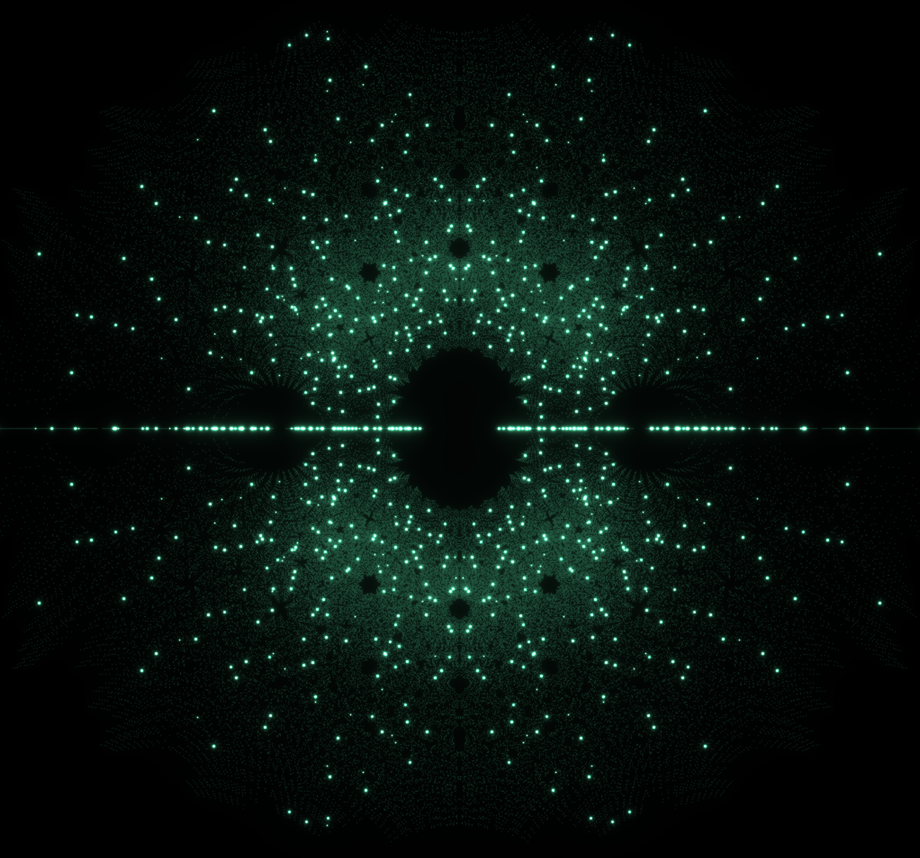}\\
  \vspace{20pt}\includegraphics[scale=.22]{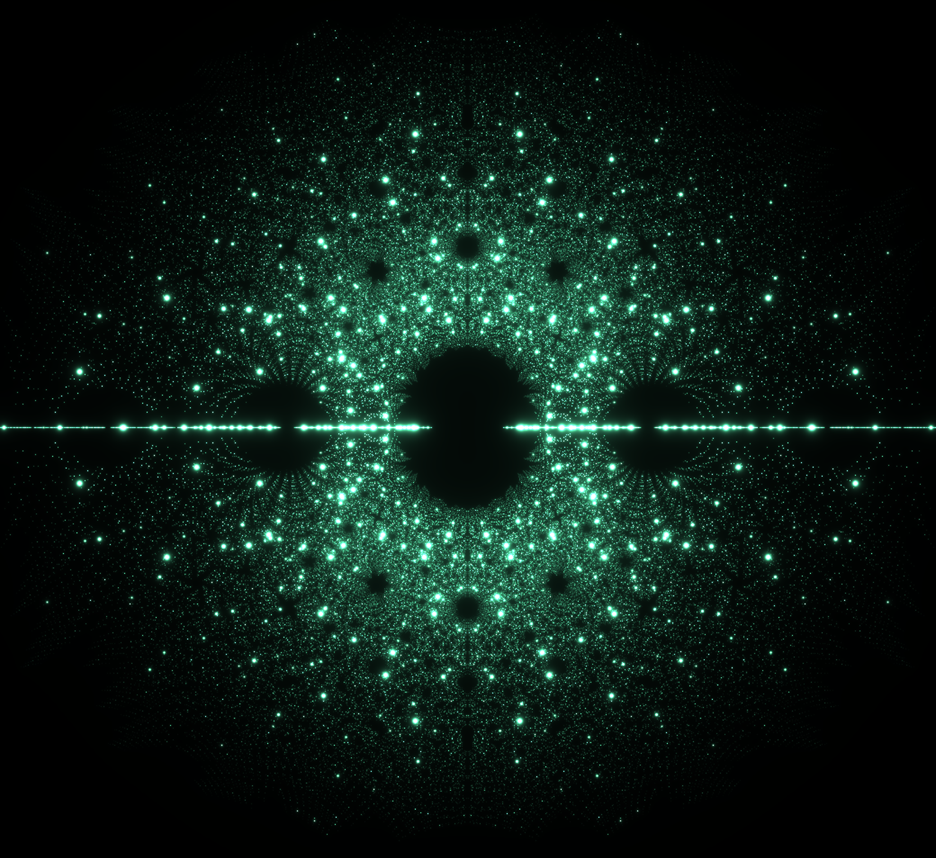}\hspace{20pt}\includegraphics[scale=.22]{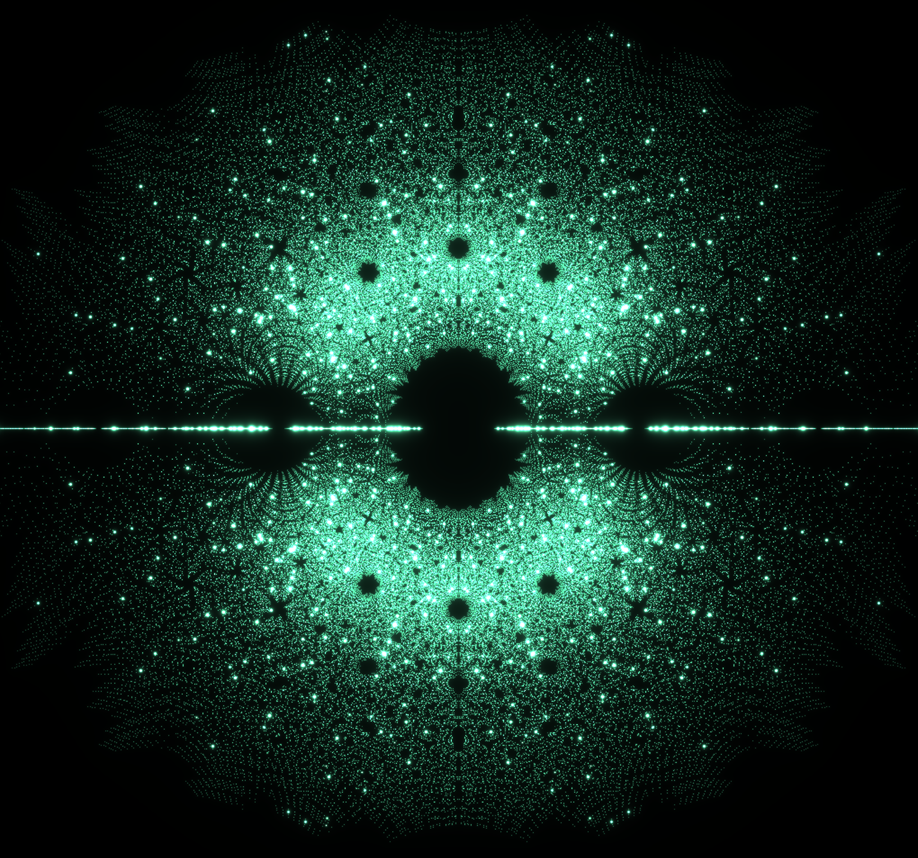}
  \caption{Algebraic integers of degree five.  \emph{Top left}: Sized by disciminant only, using $\beta_{D}$ as the brightness parameter.  \emph{Top right}: Sized by rigidity only, using $\beta_{r}^+$ as the brightness parameter.  \emph{Bottom left}: Sized by disciminant and rigidity, using $\beta_{Dr}$ as the brightness parameter.  \emph{Bottom right}: Sized by rigidity and the magnitude of the coefficient vector, using $\beta_{vr}$ as the brightness parameter.}
  \label{quintics!}
\end{figure}

We next turn our attention to quintics in Figure \ref{quintics!}.  The images are rich and beautiful, and full of interesting patterns, but other than the concentration of rigid integers at the real line, predicted by Remark \ref{RigidReality}, it isn't immediately clear how to quantitatively describe the patterns appearing when we highlight rigidity.  There does seem to be some interesting underlying structure and symmetry, and a few potentially hexagonal patterns appear where there is a denser concentration near the middle.  We hope that by expanding our coefficient ranges and amassing more data, we can fill in the space and illuminate the underlying structure more clearly.

\begin{figure}[h!]
  \centering
  \includegraphics[scale=.3]{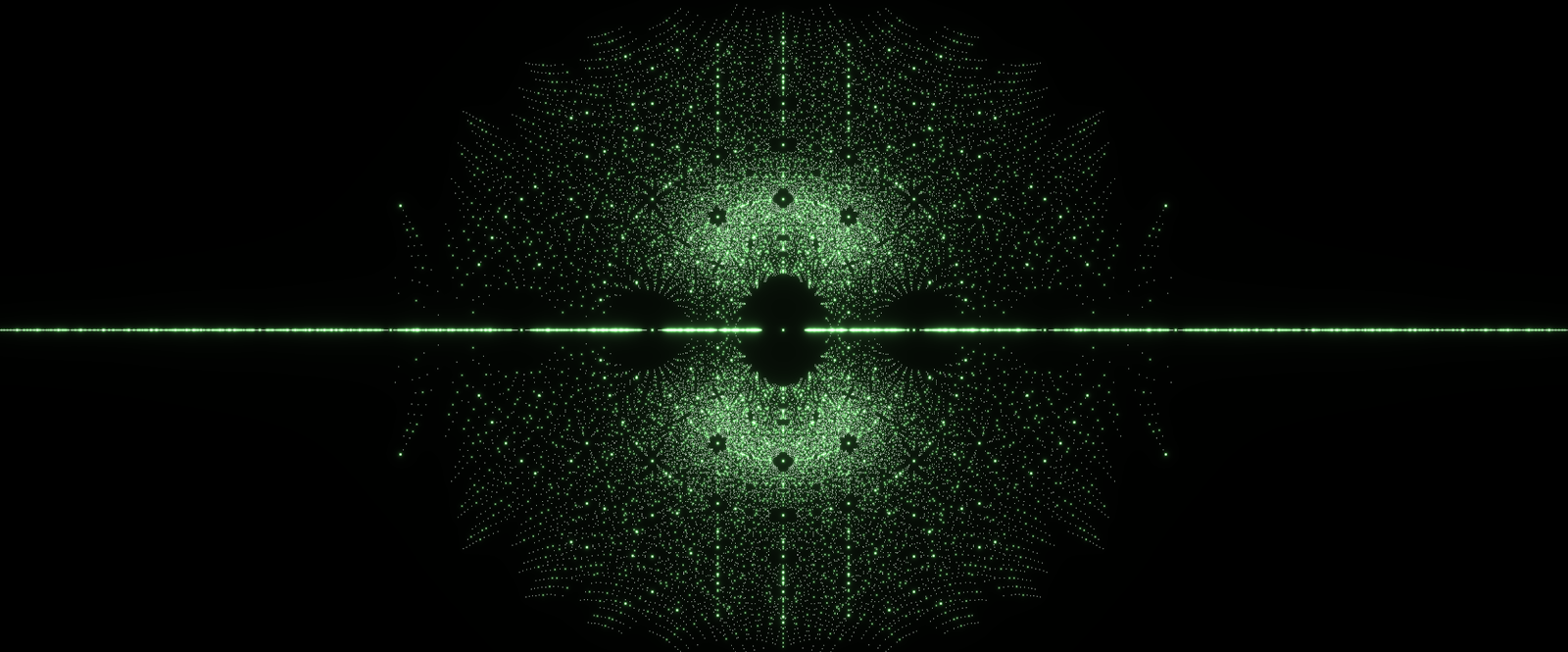}
  \caption{Algebraic integers of degrees one through four with brightness parameter $\beta(\alpha) = C(\deg(\alpha)/|G_\alpha|)^2$.}
  \label{mixedUp}
\end{figure}

To conclude we consider one more image (Figure \ref{mixedUp}) which is built from a wider coefficient range, but only including degrees one through four.  It also uses a different formula for brightness, setting $\beta(\alpha) = C(\deg(\alpha)/|G_\alpha|)^2$ (where $C$ is a constant).  We point out that since the size of Galois group is usually the degree factorial, the size of points decreases significantly as degree goes up, so the constant $C$ is adjusted for each degree.  Therefore, one should remain cautious when drawing too many conclusions about geometric relationships of rigid integers of different degrees in this image.  Nevertheless, the appearance of the elliptic orbits of rigid points about the origin which do not appear in any previous images suggest that widening our coefficient ranges should indeed produce even more intriguing results.

We invite the interested reader to use the dataset and Mathematica code available in the public github repository \cite{github} to experiment with different brightness functions and coefficient ranges.  There, one can also find a gallery of the many images we have generated by adjusting the inputs (of which we have only shared a few in this report).
\subsection{Future Work}
There are a number of angles we consider for future work stemming from this project.  The first is theoretical in nature, attempting to get complete and rigorous explanations of the patterns appearing in the generated images, extending Propositions \ref{cubics} and \ref{RatGeo}.  An analysis where algebraic integers were sized according to discriminant was studied very successfully in \cite{hst21}, and although Galois groups are somewhat more mysterious objects in general, a full explanation in the case of cubics (Proposition \ref{cubics}) and a partial explanation in the case of quartics (Proposition \ref{RatGeo}) gives a good start, leaving us optimistic about the prospects.

Another next step would be to generate more data, particularly in the quintic case, and to move on to higher degrees.  This would allow us to generate a more complete illustration of what structures rigid points assume, allowing us to formulate more precise questions.  Right now our method consists of computing the degree of the splitting field, and we were able to do it rather efficiently for the 220,000 quintic algebraic integers appearing in the quintic images we generated, but efficiency falls off quickly after that.\footnote{For reference, the computation was done on a modestly powered Dell XPS 15 7590 with an Intel Core i7-9750H (12MB Cache, up to 4.5 GHz, 6 cores).}  Going further, the process should lend itself well to parallel computing, which is an avenue we have not yet explored.  We also expect a more theoretical route could be fruitful.  Indeed, for quartics, there is a classification of Galois groups of an irreducible quartic in terms of whether the discriminant is a square, and whether the resolvent cubic is irreducible (cf. \cite[Section 14.6]{DF04}), which allows one to avoid explicitly computing the splitting field, thereby speeding up the computation substantially.  Harnessing even partial results in this direction for higher degree polynomials would allow for a much more efficient implementation.

Finally, we are interested in exploring similar projects but for different invariants associated to algebraic integers.  For example, one could consider things like ramification and inertia for various primes, or the norm or trace of the minimal polynomial, and generate different landscapes of algebraic integers with interesting geometric structures to explore.
\begin{center}
  \begin{tabular}{l c l}
    Gabriel Dorfsman-Hopkins &{ }& Shuchang Xu\\
    University of California, Berkeley && University of California, Berkeley\\
    \verb|gabrieldh@berkeley.edu| && \verb|candyxu@berkeley.edu|\\
    895 Evans Hall && \\
    Berkeley, CA, 94720 && \\
    \verb|www.gabrieldorfsmanhopkins.com| &&
  \end{tabular}
\end{center}
\bibliography{bib}{}
\bibliographystyle{alpha}
\end{document}